\documentclass[14pt]{extarticle}
\usepackage[cp1251]{inputenc}
\usepackage[russian]{babel}
\usepackage{amsmath}
\usepackage{amssymb}
\usepackage[russian]{babel} 

    \usepackage[pdftex]{graphicx}

\newcommand{\prodd}{\prod\limits}
\newcommand{\summ}{\sum\limits}
\newcommand{\D}{\mathcal{D}}
\newcommand{\F}{\mathbb{F}}
\renewcommand{\leq}{\leqslant}
\renewcommand{\geq}{\geqslant}

\begin{document}

\begin{flushleft}УДК 517
\end{flushleft}
\begin{center}
\textbf{О КВАДРАТАХ В СПЕЦИАЛЬНЫХ МНОЖЕСТВАХ КОНЕЧНОГО ПОЛЯ\footnote{Исследование выполнено за счет гранта Российского научного фонда (проект 14-11-00702).}} 
\end{center}	
\begin{center}
М.\,Р.~Габдуллин (г. Москва)
\end{center}
\begin{abstract}
Рассматривается линейное пространство, образованное элементами конечного поля $\F_q$, где $q=p^r$, над $\F_p$. Пусть $\{a_1,\ldots,a_r\}$ -- базис этого пространства. Тогда каждый элемент $x\in\F_q$ имеет единственное представление в виде $\sum_{j=1}^r c_ja_j$, где $c_j\in\F_p$. Зафиксируем множества $D_1,\ldots,D_r\subset\F_p$ и рассмотрим множество $W=W(D_1,\ldots,D_r)$ тех элементов $x\in\F_q$, для которых $c_j\in D_j$ при всех $1\leq j \leq r$. 

В работе доказана оценка на количество квадратов во множестве $W$, из которой вытекают следующие два утверждения:

1) если для некоторого $\varepsilon>0$ выполнено $\prodd_{i=1}^r|D_i| \geq (2r-1)^rp^{r(1/2+\varepsilon)}$, то справедлива асимптотическая оценка $|W\cap Q|=|W|\left(\frac12+O(p^{-\varepsilon/2})\right)$; 

2) при $\prod\limits_{i=1}^r |D_i|\geq 8(2r-1)^rp^{r/2}$ во множестве $W$ имеются ненулевые квадраты.

\bigskip

\textit{Ключевые слова:} конечные поля, квадраты, суммы характеров.   	 
\end{abstract}	

\bigskip

\section{Введение.}

При любом фиксированном $b\in\mathbb{N}$, $b\geq2$, каждое число $n\in\mathbb{N}$ единственным образом представимо в системе счисления с основанием $b$:

\begin{equation*}
n=\sum\limits_{j=0}^{r-1} c_jb^j, \,\,\,\,\,\,\, 0\leq c_j \leq b-1, \,\,\, c_{r-1}\geq 1.
\end{equation*} 
Во многих работах (см., например, \cite{1}-\cite{11.5})  изучались арифметические свойства чисел с ''пропущенными'' цифрами, т.е. тех чисел, $b$-ичная запись которых состоит из заданных цифр. 

В \cite{12} C.~Dartyge и A.~S\'ark\"ozy рассмотрели аналог этой задачи в конечных полях. Пусть $\mathbb{F}_q$ --- поле из $q=p^r$ элементов, $r\geq2$, $\{a_1,\ldots,a_r\}$ --- базис $\mathbb{F}_q$ над $\mathbb{F}_p$,  $\mathcal{D}\subset\mathbb{F}_p$. Положим
$$ W_{\D}=\left\{ x_1a_1+\ldots+x_ra_r \,|\, x_i\in \D  \right\}. 
$$
Обозначим через $Q$ множество ненулевых квадратов поля $F_q$. Положим $Q_0=Q\cup\{0\}.$ Будем считать, что $p\geq3$, так как в случае $p=2$ мы имеем $\F_q=Q_0.$

В недавней работе C.Dartyge, C.Mauduit, A.S\'ark\"ozy \cite{13} было показано, что если множество $\D$ достаточно велико, то во множестве $W_{\D}$ имеются квадраты.

\textsc{Теорема A}. \textit{Пусть $2\leq|\D|\leq p-1.$ Тогда}
$$\left||W_{\D}\cap Q_0|-\frac{|W_{\D}|}{2}\right| \leq \frac{1}{2\sqrt{q}}\left(|\D|+p\sqrt{p-|\D|}\right)^r.
$$
Эта оценка нетривиальна, если $|\D|\geq\frac{(\sqrt5-1)p}{2}(1+o_p(1))$.

В работе \cite{14} автором были доказаны следующие два утверждения, ослабляющие условие на мощность множества $\D$, гарантирующее наличие квадратов во множестве $W_{\D}$.

\textsc{Теорема B.} \textit{Пусть $2r-1\leq p^{1/2}$, $\delta=\left(\sqrt p(2r-1)\right)^{2-r}$. Тогда при \newline $|\D|\geq(1+\delta)(2r-1)p^{1/2}$ справедливо $|W_{\D}\cap Q|\geq1.$}

\bigskip

\textsc{Теорема C.} \textit{Пусть $r\geq20,$  $C(r)=\exp\left(\frac{4\log r+8}{r}\right)=1+o(1), \, r\to\infty$. Тогда при $|\D|\geq  C(r) p^{\frac12}\exp\left(\frac{\log p+4\log\log p}{r}\right)$ справедливо $|W_{\D}\cap Q|\geq1$.}

\bigskip

В частности, из теоремы C следует, что при большом $r$ во множестве $W_{\D}$ есть квадраты уже при $|\D|>p^{1/2}$. Отметим, что при $r \gg \frac{\log p}{\log\log p}$, более точный результат дает теорема C, а иначе --- теорема B.

В работе R.Dietmann, C.Elsholtz, I.E.Shaprlinski \cite{15}, была рассмотрена более общая задача. Пусть $D_1,\ldots,D_r$ -- подмножества  $\mathbb{F}_p$. Положим
$$ W=W(D_1,\ldots,D_r)=\left\{ x_1a_1+\ldots+x_ra_r \,|\, x_i\in D_i  \right\}. 
$$
Авторы работы \cite{15} отмечают, что доказательство теоремы А \cite{13} переносится на случай, когда множества $D_i$ различны, а именно, при $\min\limits_{1\leq i\leq r} |D_i| \geq \frac{(\sqrt5-1)p}{2}(1+o_p(1))$ справедливо $|W\cap Q_0|\geq 1$, и доказывают более сильное утверждение.  

\textsc{Теорема D} (\cite{15}, теорема 3.5). \textit{Для любого $\varepsilon>0$ существует $\delta>0$ такое, что для любых множеств $D_1,\ldots,D_r$, удовлетворяющих условиям}
$$ \prod_{i=1}^r |D_i| \geq p^{(1/2+\varepsilon)r^2/(r-1)}
$$
\textit{и}
$$ \min\limits_{1\leq i\leq r} |D_i| \geq p^{\varepsilon}
$$
\textit{справедливо} $|W\cap Q_0|=\left(\frac12+O(p^{-\delta})\right)|W|$.

По аналогии с работой \cite{15}, теорема B также может быть перенесена на случай различных множеств $D_i$. В настоящей работе будет доказана следующая 

\textbf{Теорема.} \textit{Справедлива оценка}
\begin{multline}
\left| |W\cap Q|-\frac{|W|}{2} \right| \leq \\
\frac12\left(|W|^{1-1/(2r)}p^{1/4}(2r-1)^{1/2} + |W|^{1/(2r)}(\frac14p^{3/4}r^{3/2}+p^{1/2})+1\right). \label{T}
\end{multline}

Из этой теоремы вытекает аналог теоремы D, а также теорема о достаточных условиях существования квадратов во множестве $W$.
   
\textsc{Следствие 1.} \textit{Пусть для некоторого $\varepsilon>0$ выполнено}
$$ \prodd_{i=1}^r|D_i|\geq (2r-1)^rp^{r(1/2+\varepsilon)}.
$$ 

\textit{Тогда} $|W\cap Q|=|W|\left(\frac12+O(p^{-\varepsilon/2})\right)$, \textit{причем постоянная в знаке $O$ абсолютна.}. 

\smallskip
Отметим, что следствие 1 усиливает теорему D при фиксированном $r$ (так как в нём отсутствует требование $\min\limits_{1\leq i\leq r} |D_i| \geq p^{\varepsilon}$). 
 
\smallskip 
\textsc{Следствие 2.} \textit{Пусть $\prodd_{i=1}^r|D_i|\geq 8(2r-1)^rp^{r/2}$. Тогда $|W\cap Q|\geq 1.$}

\bigskip

Доказательство теоремы будет изложено в разделе 2; оно основывается на оценке сумм характеров специального вида, полученной D.Wan в \cite{15.5} и сформулированной в удобном для нас виде A.~Winterhof в \cite{16}.

\textsc{Лемма}. \textit{Пусть $\chi$ -- мультипликативный характер порядка $s$ в $\F_{q}$ и $\alpha,\beta\in\F_{q}$ -- несопряжённые порождающие элементы $\F_{q}$ над $\F_p$. Тогда}

$$\left|\sum\limits_{\xi\in\F_p}\chi\left((\xi+\alpha)(\xi+\beta)^{s-1}\right)\right| \leq (2r-1)p^{1/2}.
$$

\bigskip

Следствия 1 и 2 будут доказаны в разделе 3. 

Автор признателен С.\,В.~Конягину и рецензенту за полезные обсуждения результатов.

\section{Доказательство теоремы.} Через $\chi$ обозначим квадратичный характер на $\mathbb{F}_q$; cчитаем, что $\chi(0)=0$. Пусть $0$ не принадлежит некоторому $D_i$. Тогда

$$|W\cap Q|=\frac12\sum\limits_{x\in W} (1+\chi(x))=\frac12|W|+\frac12\sum\limits_{x\in W} \chi(x)=\frac12|W|+\frac12\sum\limits_{x\in W} \chi(x).
$$
Если же $0\in D_i$ при всех $i$, то
$$|W\cap Q|=\frac12\sum\limits_{x\in W\setminus\{0\}} (1+\chi(x))=\frac12(|W|-1)+\frac12\sum\limits_{x\in W} \chi(x).
$$
 Таким образом, всегда справедливо неравенство
\begin{equation}
\left| |W\cap Q|-\frac{|W|}{2} \right| \leq \frac12+\frac12\left|\summ_{x\in W} \chi(x)\right|,  \label{0}
\end{equation} 
и для доказательства теоремы нужно оценить сумму характеров $\left|\sum\limits_{x\in W}\chi(x)\right|$. Положим $D=D_2\times\ldots\times D_r$ и $b_j=a_j/a_1$. Тогда $b_1=1$ и $\{1,b_2,\ldots,b_r\}$ --- базис. Имеем
\begin{equation}
\left| \sum\limits_{x\in W} \chi(x)\right| \leq \sum\limits_{c_1\in D_1}\left| \sum\limits_{c_j\in D_j,\, j\geq2} \chi(c_1a_1+\ldots+c_ra_r)\right| \leq |D_1|^{1/2}A^{1/2}, \label{1}
\end{equation}
где
$$ A=
\sum\limits_{c_1\in D_1}\left| \sum\limits_{(c_2,\ldots,c_r)\in D} \chi(c_1+c_2b_2\ldots+c_rb_r)\right|^2.
$$
Пусть $L_d$ --- множество тех наборов $(c_2,\ldots,c_r)\in D$, для которых элемент $c_2b_2+\ldots+c_rb_r$ лежит в подполе порядка $p^d$ и не лежит ни в каком подполе меньшего порядка. Ясно, что $D=\bigsqcup\limits_{d|r} L_d$, причем $L_1=\{0\}$, если $0\in D_i$ при всех $2\leq i \leq r$, и $L_1=\emptyset$ иначе. Для $d|r$ определим функцию $f_d(x) \colon D_1\to\mathbb{C}$,  $f_d(c)=\sum\limits_{(c_2,\ldots,c_r)\in L_d} \chi(c+c_2b_2+\ldots+c_rb_r)$. Напомним, что $l_2$-норма функции $g \colon D_1\to\mathbb{C}$ определяется как $\| g \|_2 = \left(\sum\limits_{x\in D_1}|g(x)|^2 \right)^{1/2}$. Тогда в силу неравенства треугольника

\begin{equation}A^{1/2}=\left\| \sum_{d|r} f_d \right\|_2 \leq \sum\limits_{d|r} \|f_d\|_2 = \sum\limits_{d|r} A_d^{1/2}. \label{A^{1/2}}
\end{equation}
где
$$ A_d= \sum\limits_{x\in D_1}\left| \sum\limits_{(c_2,\ldots,c_r)\in L_d} \chi(x+c_2b_2+\ldots+c_rb_r)\right|^2 .
$$
По определению множества $L_d$ при любом наборе $(c_2,\ldots,c_r)\in L_d$ элемент $c_2b_2+\ldots+c_rb_r$ порождает подполе порядка $p^d$. Учитывая, что каждый такой элемент имеет не более $d$ сопряженных, и применяя лемму к парам несопряженных элементов, при $d>1$ имеем
\begin{multline*} A_d\leq \sum\limits_{(c_2,\ldots,c_r), (c'_2,\ldots,c'_r)\in L_d} \left|\sum\limits_{x\in F_p}  \chi(x+c_2b_2+\ldots+c_rb_r)\overline{\chi}(x+c'_2b_2+\ldots+c'_rb_r)\right|
\leq \\
\sum\limits_{(c_2,\ldots,c_r)\in L_d} \left(dp+(|L_d|-d)(2d-1)p^{1/2}\right)\leq(2d-1)p^{1/2}|L_d|^2+dp|L_d|.
\end{multline*}
Кроме того, $A_1\leq |D_1|\leq p.$ Обозначим $\mathcal{J}=\{d|r : d>1 \,\mbox{и}\, L_d \neq\emptyset \}$. Тогда при $d\in\mathcal{J}$ в силу неравенства $\sqrt{A+B}\leq \sqrt{A}\left(1+\frac{B}{2A}\right)$, верного при всех положительных $A$ и $B$, получаем
$$ A_d^{1/2}\leq (2d-1)^{1/2}p^{1/4}|L_d|+\frac{dp^{3/4}}{2(2d-1)^{1/2}}.
$$
Из этой оценки и неравенства (\ref{A^{1/2}}) имеем
$$ A^{1/2}\leq p^{1/4}S_1+\frac12p^{3/4}S_2+p^{1/2},
$$
где
\begin{gather*}
S_1=\sum\limits_{d\in \mathcal{J}} (2d-1)^{1/2}|L_d|,\quad
S_2=\sum\limits_{d\in\mathcal{J}}\frac{d}{(2d-1)^{1/2}}.
\end{gather*}
Учитывая, что $\sum\limits_{d|r} |L_d|=|D_2|\ldots|D_r|$, получаем
$$
S_1\leq(2r-1)^{1/2}|D_2|\ldots|D_r|,\quad
S_2\leq\sum\limits_{d\in \mathcal{J}} d^{1/2} \leq \frac12r^{3/2}.
$$
(Последняя оценка проверяется непосредственно при $2\leq r\leq7$, а при $r \geq 8$ вытекает из неравенств $r^{1/2}\leq\frac16r^{3/2}$ и $\sum\limits_{d\leq r/2} d^{1/2}\leq \frac23(r/2+1)^{3/2}\leq \frac13r^{3/2}.$)
Значит,
$$A^{1/2}\leq p^{1/4}(2r-1)^{1/2}|D_2|\ldots|D_r|+\frac14p^{3/4}r^{3/2}+p^{1/2}.
$$
Подставляя последнее неравенство в (\ref{1}), получим
\begin{multline*}
\left| \sum\limits_{x\in W} \chi(x)\right| \leq |D_1|^{1/2}\left(p^{1/4}(2r-1)^{1/2}|D_2|\ldots|D_r| + \frac14p^{3/4}r^{3/2}+p^{1/2}\right) =\\ 
|D_1|^{-1/2}p^{1/4}(2r-1)^{1/2}|W|+|D_1|^{1/2}\left(\frac14p^{3/4}r^{3/2}+p^{1/2}\right)
\end{multline*}
Аналогично получается оценка
\begin{multline*}
\left| \sum\limits_{x\in W} \chi(x)\right| \leq
|D_i|^{-1/2}p^{1/4}(2r-1)^{1/2}|W|+|D_i|^{1/2}\left(\frac14p^{3/4}r^{3/2}+p^{1/2}\right)
\end{multline*}
Выберем $i$ так, чтобы эта оценка была наилучшей. Рассмотрим функцию \\$H(x)=x^{-1/2}p^{1/4}(2r-1)^{1/2}|W|+x^{1/2}\left(\frac14p^{3/4}r^{3/2}+p^{1/2}\right)$, $1\leq x \leq |W|$. Имеем 

$$H'(x)=\frac{\frac14p^{3/4}r^{3/2}+p^{1/2}}{2x^{1/2}} - \frac{p^{1/4}(2r-1)^{1/2}|W|}{2x^{3/2}},
$$
и $H'(x)=0$ при $x=4|W|p^{-1/2}\frac{(2r-1)^{1/2}}{r^{3/2}+4p^{-1/4}}\geq 2|W|p^{-1/2}r^{-1}\geq |W|^{1/2}$ при $|W|\geq 8p^{r/2}(2r-1)^r$. Таким образом, минимум функции $H(x)$ достигается при $x\geq|W|^{1/2}\geq|W|^{1/r}$. Мы можем выбрать $i$ так, что $|D_i|\geq|W|^{1/r}$; так как $H'(x)<0$ при $|W|^{1/r}\leq x \leq |W|^{1/2}$, то можно гарантировать оценку
\begin{multline*} \left| \sum\limits_{x\in W} \chi(x)\right| \leq H(|W|^{1/r})= \\ |W|^{1-1/(2r)}p^{1/4}(2r-1)^{1/2}+|W|^{1/(2r)}\left(\frac14p^{3/4}r^{3/2}+p^{1/2}\right).
\end{multline*}
Из неравенства (\ref{0}) и последней оценки вытекает утверждение теоремы. 

\section{Доказательство следствий.}

\textbf{Доказательство следствия 1.} Из теоремы следует, что

$$|W\cap Q|=|W|\left(\frac12+O\left( |W|^{-1/(2r)}p^{1/4}(2r-1)^{1/2}+|W|^{-(2r-1)/(2r)}p^{3/4}r^{3/2} \right)\right).
$$
При $|W|\geq (2r-1)^rp^{r(1/2+\varepsilon)}$ имеем $|W|^{-1/(2r)}\leq (2r-1)^{-1/2}p^{-1/4}p^{-\varepsilon/2}$ и
$$|W\cap Q|=|W|\left(\frac12+O\left(p^{-\varepsilon/2}+r^{2-r}p^{1-r/2}p^{-(2r-1)\varepsilon/2}\right)\right).
$$
Так как $r\geq2$, то отсюда вытекает утверждение следствия 1.

\textbf{Доказательство следствия 2.} Во множестве $W$ есть квадраты, если правая часть неравенства (\ref{T}) строго меньше, чем $|W|-1$. Это равносильно условию
$$|W|^{1-1/r}\left(|W|^{1/(2r)}-p^{1/4}(2r-1)^{1/2}\right) > \frac14p^{3/4}r^{3/2}+p^{1/2}+|W|^{-1/(2r)}.
$$
Покажем, что последнее неравенство выполнено при $|W|\geq 8(2r-1)^rp^{r/2}$. Имеем $|W|^{1/(2r)}\geq \left(e^2(2r-1)^rp^{r/2}\right)^{1/(2r)}\geq (1+1/r)(2r-1)^{1/2}p^{1/4}$ и
\begin{multline*}
|W|^{1-1/r}\left(|W|^{1/(2r)}-(2r-1)^{1/2}p^{1/4}\right)\geq\frac1r(2r-1)^{1/2}p^{1/4}(2r-1)^{r-1}p^{(r-1)/2} =\\
\frac1r(2r-1)^{r-1/2}p^{(r-1)/2+1/4} \geq \frac34p^{3/4}r^{3/2}>\frac14p^{3/4}r^{3/2}+p^{1/2}+|W|^{-1/(2r)}.
\end{multline*}
Здесь предпоследнее неравенство очевидно при $r\geq 3$ и легко проверяется при $r=2$; последнее неравенство следует из неравенства $\frac12p^{3/4}r^{3/2}> p^{1/2}+(2r-1)^{-1/2}p^{-1/4}$, верного при всех $p\geq3$, $r\geq2$. 

Следствие доказано.

\section{Заключение.}  В работе доказана оценка на количество квадратов во множестве $W$, дающая аналог результата из работы \cite{15}, а также достаточные условия на существование квадратов во множестве $W$. 

Этот результат является обобщением теоремы B предыдущей работы автора \cite{14}.

\flushleft{\textbf{М.\,Р.~Габдуллин} \\
Московский государственный университет им. \\
Ломоносова\\
Институт математики и механики УрО РАН. \\
\textit{E-mail:} Gabdullin.Mikhail@ya.ru}

\end{document}